\documentclass{amsart}
\usepackage{graphicx}

\newtheorem{theorem}{Theorem}[section]
\newtheorem{lemma}[theorem]{Lemma}

\theoremstyle{definition}

\newtheorem{example}[theorem]{Example}

\theoremstyle{remark}

\numberwithin{equation}{section}

\begin{document}

\title{Extension of Lyapunov's Convexity Theorem to Subranges}

\author{Peng Dai}
\address{Department of Applied Mathematics and Statistics, Stony Brook University, Stony Brook, NY 11794-3600}
\email{Peng.Dai@stonybrook.edu}

\author{Eugene A. Feinberg}
\address{Department of Applied Mathematics and Statistics, Stony Brook University, Stony Brook, NY 11794-3600}
\email{Eugene.Feinberg@sunysb.edu}

\thanks{This research was partially supported by NSF grants CMMI-0900206 and CMMI-0928490.}



\date{\today}


\begin{abstract}
Consider a measurable space with a finite vector measure. This
measure defines a mapping of the $\sigma$-field into a Euclidean
space. According to Lyapunov's convexity theorem, the range of
this mapping is compact and, if the measure is atomless, this
range is convex. Similar ranges are also defined for measurable
subsets of the space.  We show that the union of the ranges of all
subsets having the same given vector measure is also compact and,
if the measure is atomless, it is convex. We further provide a
geometrically constructed convex compactum in the Euclidean space
that contains this union. The equality of these two sets, that
holds for two-dimensional measures, can be violated in higher
dimensions.
\end{abstract}
\maketitle

\section{Introduction}\label{sec:intro}

Let $\left({X},\mathcal{F}\right)$ be a measurable space and
${\mu}=\left(\mu_{1},...,\mu_{m}\right)$, $m=1,2,\ldots,$ be a
finite vector measure on it. For each ${Y}\in \mathcal{F}$
consider the range ${R}_{{\mu}}\left({Y}\right)=\left\{
{\mu}\left({Z}\right):{Z}\in\mathcal{F},{Z}\subset{Y}\right\}\subset
\mathbb{R}^m$ of the vector measures of all its measurable subsets
$Z$. Lyapunov's convexity theorem~\cite{Liapounoff:1940} states
that the range ${R}_{{\mu}}\left({X}\right)$ is compact and
furthermore, if $\mu$ is atomless, this range is convex. Of
course, this is also true for any $Y \in \mathcal{F}$. We recall
that a measure $\nu$ is called atomless if for each
${Z}\in\mathcal{F}$ such that $\nu\left({Z}\right)>0$, there
exists ${Z}'\in\mathcal{F}$ such that  ${Z}'\subset{Z}$ and
$0<\nu\left({Z}'\right)<\nu\left({Z}\right)$. A vector measure
${\mu}=\left(\mu_{1},...,\mu_{m}\right)$, is called atomless if
each measure $\mu_{i}$, $i=1\dots m$, is atomless. For a review of
Lyapunov's convexity theorem and its applications see
\cite{Olech}.

 Let $\mathcal{S}^{{p}}_\mu(X)$
 be the set of all measurable subsets of ${X}$ with the vector measure ${p}\in R_\mu(X)$,
\[
\mathcal{S}^{{p}}_\mu\left({X}\right)=\left\{ {Y}\in\mathcal{F}:
{\mu}\left({Y}\right)={p}\right\} .
\]
Of course, $\mathcal{S}^{{p}}_\mu\left({X}\right)=\emptyset$, if
$p\notin R_\mu(X)$.
 For ${p}\in \mathbb{R}^m$  consider the union of
the ranges of all subsets of ${X}$ with the vector
measure ${p}$,
\[
{R}^{{p}}_{\mu} \left({X}\right)=\bigcup \limits
_{{Y}\in\mathcal{S}^{{p}}_\mu\left({X}\right)}{R}_\mu\left({Y}\right).
\]
In particular, ${R}^{{p}}_{\mu} \left({X}\right)=\emptyset, $ if
$p\notin R_\mu(X)$, and ${R}^{{\mu(X)}}_{\mu}
\left({X}\right)={R}_{\mu} \left({X}\right). $

 Since the relation $Y_1=Y_2$ ($\mu$-everywhere) is an equivalence
relation on $\mathcal F$, it partitions any subset of $\mathcal F$
into equivalence classes.
 For an atomless
$\mu$, Lyapunov \cite[Theorem III]{Liapounoff:1940} proved that:
(i) $\mathcal{S}^{{p}}_\mu\left({X}\right)$  consists of one
equivalence class if and only if $p$ is an extreme point of
$R_\mu(X)$, and (ii) if $p\in R_\mu(X)$ is not an extreme point of
$R_\mu^p(X)$, then the set of equivalence classes in
$\mathcal{S}^{{p}}_\mu\left({X}\right)$ has cardinality of the
continuum.

In general, a union of an infinite number of compact convex sets
may be neither closed nor convex.  As follows from Dai and
Feinberg \cite{Dai:2010}, the set $R_\mu^p(X)$ is a convex
compactum, if $m=2$ and $\mu$ is atomless.  This fact follows from
stronger results that hold for $m=2.$   For $m=2$ and atomless
$\mu$,  Dai and Feinberg \cite[Theorem~2.3]{Dai:2010} showed that
there exists a set $Z^*\in \mathcal{S}^{{p}}_\mu\left({X}\right)$,
called a maximal subset,   such that
\begin{equation} \label{eq:2D}
R_\mu(Z^*)=R^p_\mu(X),
\end{equation}
and, in addition, the following equality holds
\begin{equation} \label{eq:2DA}
R^p_\mu(X)=Q^p_\mu(X),
\end{equation}
where $Q^p_\mu(X)$ is the intersection of $R_\mu(X)$ with its shift by a vector $- ({ \mu(X) -p})$,
\begin{equation}\label{eq:defQ}
Q^{{p}}_{\mu}\left({X}\right)=\left({R}_{{\mu}}\left({X}\right)-\left\{
{\mu}\left({X}\right)-{p}\right\}
\right)\cap{R}_{{\mu}}\left({X}\right),
\end{equation}
with ${S}_1-{S}_2=\{ {q}-{r}:{q}\in {S}_1,{r}\in{S}_2\}$ for
${S}_1,{S}_2\subset \mathbb{R}^m$. In particular,
${R}_{{\mu}}\left({X}\right)-\{{r}\}$ is a parallel shift of
${R}_{{\mu}}\left({X}\right)$ by $-{r}\in \mathbb{R}^m$. Examples
\ref{ex3} and \ref{ex4} below demonstrate that  equalities
(\ref{eq:2D}) and (\ref{eq:2DA}) may not hold when $\mu$ is not
atomless even for $m=1.$

Each of equalities (\ref{eq:2D}) and (\ref{eq:2DA}) implies that
$R^p_\mu(X)$ is convex and compact.  However, \cite[Example
4.2]{Dai:2010} demonstrates that a maximal set $Z^*$ may not exist
for an atomless vector measure $\mu$ when $m>2.$

In this paper, we prove (Theorem 1) that for any natural number
$m$ the set $R_\mu^p(X)$ is compact and, if $\mu$ is atomless,
this set is convex. This is a generalization of Lyapunov's
convexity theorem, which is a particular case of this statement
for  $p=\mu(X).$ We also prove that  $R_{\mu}^{p}(X)\subset
Q_{\mu}^{p}(X)$ (Theorem 2).  Example~\ref{ex:RinQ} demonstrates
that it is possible that
 equality  (\ref{eq:2DA}) may not hold  when $m>2$ and $\mu$ is
 atomless.

Lyapunov's convexity theorem is relevant to purification of
transition probabilities discovered by Dvoretzky, Wald, and
Wolfowitz
\cite{Dvoretzky:1950, Dvoretzky:1951}. 
Let $(A,{\mathcal{A}})$ be a measurable space and $\pi$ be a
transition probability from $(X,{\mathcal F})$ to $A$; that is,
$\pi(B|x)$ is a measurable function on $(X,{\mathcal{F}})$ for any
$B\in \mathcal{A}$ and $\pi(\cdot|x)$ is a probability measure on
$(A,{\mathcal{A}})$ for any $x\in X$. According to Dvoretzky,
Wald, and Wolfowitz \cite{Dvoretzky:1950, Dvoretzky:1951}, two
transition probabilities $\pi_1$ and $\pi_2$ are called strongly
equivalent if
\begin{equation}\label{eq:SETEQ}
\int_{{X}}\pi_{1}\left(B|x\right)\mu_i\left(dx\right)=
\int_{{X}}\pi_{2}\left(B|x\right)\mu_i\left(dx\right),
\qquad i=1,\dots,m, \quad B\in {\mathcal{A}}.\end{equation}

A transition probability $\pi$ is called pure if each measure
$\pi(\cdot|x)$ is concentrated at one point.  A pure transition
probability $\pi$ is defined by a measurable mapping $\varphi:X\to
A$ such that $\pi(B|x)=I\{\varphi(x)\in B\}$ for all $B\in
\mathcal{A}.$  According to the contemporary terminology, a
transition probability can be purified if  it is strongly
equivalent to a pure transition probability.

For a finite set $A$, Dvoretzky, Wald, and Wolfowitz
\cite{Dvoretzky:1950,Dvoretzky:1951} proved that any transition
probability can be purified, if the measure $\mu$ is atomless.
Edwards \cite[Theorem 4.5]{Edwards:1987} generalized this result
to the case of a countable set $A$. Khan and Rath \cite[Theorem
2]{Khan:2009} gave another proof of this generalization. Loeb and
Sun \cite[Example 2.7]{Loeb:2006} constructed an elegant example
when a transition probability cannot be purified for $m=2$,
$X=[0,1]$,  $A=[-1,1]$, and atomless $\mu$. However, purification
holds for a countable set of atomless, finite, signed Loeb
measures, when $A$ is a complete separable metric space
\cite[Corollary~2.6]{Loeb:2006}.  Podczeck \cite{Podczeck:2009}
proved that purification holds for a countable set of finite
signed measures $\mu_k$ absolutely continuous with respect to a
measure $\mu$, when $(X,{\mathcal F},\mu)$ is a super-atomless
probability space and\ $A$ is a compact metric space.

We also mention that for a finite set $A$, atomless measure $\mu$,
and measurable nonnegative functions $f_i$, $i=1,\ldots,m,$ on
$X\times A,$ Dvoterzky, Wald, and Walfowitz \cite{Dvoretzky:1951,
Dvoretzky:1951a}  proved that for any transition probability there
exists an equivalent pure transition probability. Feinberg and
Piunovskiy \cite{Feinberg:2006} proved that this is true for
standard Borel spaces $X$ and $A$.  We recall that two transition
probabilities $\pi_1$ and $\pi_2$ are called equivalent
\cite{Dvoretzky:1951, Dvoretzky:1951a} if
\begin{equation*}
\int_{{X}}\int_{{A}}f_i(x,a)\pi_{1}\left(da|x\right)\mu_i\left(dx\right)=
\int_{{X}}\int_{{A}}f_i(x,a)\pi_{2}\left(da|x\right)\mu_i\left(dx\right),
\qquad i=1,\dots,m.\end{equation*}

For a countable set $A$ and atomless $\mu$, define vectors
\begin{equation} \label{e:trpr}
              p^a=\int_X\pi(a|x)\mu(dx),\qquad\qquad a\in A.
\end{equation}
Since purification holds for an atomless $\mu$ and a countable $A$
\cite{Edwards:1987}, 
vectors  $p^a$, $a\in A$, can be presented as in (\ref{e:trpr}),
if and only if there exists a partition $\{X^a:\, X^a\in {\mathcal
F}, a\in A\}$ of the set $X$ such that $\mu(X^a)=p^a$ for all
$a\in A.$ In fact, in this form the purification theorem was
presented by Dvoretzky, Wald, and Wolfowitz \cite{Dvoretzky:1951,
Dvoretzky:1951a}   for a finite set $A$.

For $m=2$, an atomless $\mu$, and a countable $A$,  Dai and
Feinberg \cite[Theorem 2.5]{Dai:2010},  provided a necessary and
sufficient condition that for a set of $m$-dimensional vectors
$\{p^a:\, a\in A\}$  there exits the above described partition.
This condition is that:
\begin{equation}\label{eq:necsuf}
{\rm (i)}\  \sum_{a \in
A} {p^a} = \mu(X),\  {\rm and\  (ii)}\  \sum_{a \in B} {p^a} \in
R_\mu (X)\ {\rm for\  any\  finite\  subset\ } B \subset
A.
\end{equation}
Obviously, (\ref{eq:necsuf}) is a necessary condition for any
natural number $m$. Example~\ref{ex:NoPart} below implies that
this condition is not sufficient for an atomless $\mu$, when $m>2$
and $A$ consists of more than two  points.

\section{Main results}\label{sec:results}

\begin{theorem}\label{thm:cc}
For any vector ${p}\in{R}_{{\mu}}\left({X}\right)$, the set
$R_{\mu}^{p}(X)$ is  compact and, in addition, if the vector
measure $\mu$ is atomless, this set  is convex.
\end{theorem}
\begin{proof}
We say that a partition is measurable, if all its elements are
measurable sets. Consider the set
\[
V_{\mu,3}(X)=\left\{
\left(\mu(S_{1}),\mu(S_{2}),\mu(S_{3})\right):\,  \left\{
S_{1},S_{2},S_{3}\right\} \textrm{ is a measurable partition of
}X\right\}.
\] According to
Dvoretzky, Wald, and Wolfowitz \cite[Theorems~1 and
4]{Dvoretzky:1951a}, $V_{\mu,3}(X)$ is  compact and, if $\mu$ is
atomless, this set is convex. Now let
\[
W_{\mu}^{p}(X)=\left\{
\left(s_{1},s_{2},s_{3}\right):\left(s_{1},s_{2},s_{3}\right)\in
V_{\mu,3}(X), \; s_{3}=\mu(X)-p, \; s_{1}+s_{2}=p\right\}.
\]
This set is compact and, if $\mu$ is atomless, it is convex.  This
is true, because $W_{\mu}^{p}(X)$ is an intersection of
$V_{\mu,3}(X)$ and two planes in $\mathbb{R}^{3m}$.  These planes
are defined by the equations $s_{3}=\mu(X)-p$ and  $
s_{1}+s_{2}=p$ respectively. We further define
\[
S_{\mu}^{p}(X)=\left\{ s_{1}:\left(s_{1},s_{2},s_{3}\right)\in W_{\mu}^{p}(X)\right\}.
\]
Since $S_{\mu}^{p}(X)$ is a projection of $W_{\mu}^{p}(X)$,
the set $S_{\mu}^{p}(X)$  is  compact and, if $\mu$ is
atomless, it is convex.

The last step of the proof is to show that
$S_{\mu}^{p}(X)=R_{\mu}^{p}(X)$ by establishing that (i)
$S_{\mu}^{p}(X)\subset R_{\mu}^{p}(X)$, and (ii)
$S_{\mu}^{p}(X)\supset R_{\mu}^{p}(X)$. Indeed, for (i), for any
$s_{1}\in S_{\mu}^{p}(X)$, there exists $(s_{1},s_{2},s_{3})\in
W_{\mu}^{p}(X)$ or equivalently there exists a measurable
partition $\left\{ S_{1},S_{2},S_{3}\right\} $ of $X$ such that
$\mu(S_{3})=\mu(X)-p$ and $\mu(S_{1})+\mu(S_{2})=p$. Let
$Z=S_{1}\cup S_{2}$. Then $\mu(Z)=p$, $s_{1}\in R_{\mu}(Z)$, and
thus $s_{1}\in R_{\mu}^{p}(X)$. For (ii), for any $s_{1}\in
R_{\mu}^{p}(X)$, there exists a set $Z\in\mathcal{F}$, such that
$\mu(Z)=p$ and $s_{1}\in R_{\mu}(Z)$, which further implies that
there exists a measurable subset $S_{1}$ of $Z$ such that
$\mu(S_{1})=s_{1}$. Let $S_{2}=Z\setminus S_{1}$ and
$S_{3}=X\setminus Z$. Then $\mu(S_{1})+\mu(S_{2})=p$ and
$\mu(S_{3})=\mu(X)-p$, which further implies that
$\left(s_{1},\mu(S_{2}),\mu(S_{3})\right)\in W_{\mu}^{p}(X)$. Thus
$s_{1}\in S_{\mu}^{p}(X)$.
%
\end{proof}

\begin{theorem}\label{thm:enclosure}
$R_{\mu}^{p}(X)\subset Q_{\mu}^{p}(X)$ for any vector
${p}\in{R}_{{\mu}}\left({X}\right)$.
\end{theorem}
Recall that $R_{\mu}^{p}(X) = Q_{\mu}^{p}(X)$ when $m=2$ and $\mu$
is atomless; see (\ref{eq:2DA}). The proof of
Theorem~\ref{thm:enclosure} uses the following lemma.

%

\begin{lemma}\label{lem:censymsub} (\cite[Lemma 3.3]{Dai:2010})
For any vector $p\in R(X),$ each of the sets
${R}_\mu^{{p}}\left({X}\right)$ and $Q_\mu^{{p}}\left({X}\right)$
is centrally symmetric with the center $\frac{1}{2}{p}$.
\end{lemma}

  Though it is assumed in \cite{Dai:2010}
that the measure $\mu$ is atomless, this assumption is not used in
the proofs of Lemmas 3.1-3.3 in \cite{Dai:2010}.

%

\begin{proof}[Proof of Theorem~\ref{thm:enclosure}]
Let ${q}\in{R}^{{p}}_\mu\left({X}\right)$. Since $R_\mu ^p (X)
\subset R_\mu (X)$, then $q \in R_\mu (X)$. Furthermore, in view
of Lemma~\ref{lem:censymsub}, $p-q \in
{R}^{{p}}_\mu\left({X}\right)$. Therefore, $p-q\in R_\mu(X).$
Since  $R_\mu(X)$ is centrally symmetric with the center
$\frac{1}{2}\mu(X),$ then $R_\mu (X) = \left\{\mu(X) \right\} -
R_\mu(X)$ and $p-q \in R_\mu (X) = \left\{\mu(X) \right\} -
R_\mu(X)$. Therefore, ${q}\in{R}_{{\mu}}\left({X}\right)-\left\{
{\mu}\left({X}\right)-{p}\right\} $. As follows from  the
definition of $Q^{{p}}_{\mu}\left({X}\right)$ in (\ref{eq:defQ}),
${q}\in Q^{{p}}_{\mu}\left({X}\right)$.
\end{proof}

\section{Counterexamples}\label{sec:counter}

The first counterexample shows that  equality  (\ref{eq:2DA}) may
not hold when $\mu$ is atomless and $m>2$. In particular, the
inclusion in Theorem \ref{thm:enclosure} cannot be substituted
with the equality.

\begin{example}\label{ex:RinQ}

Consider the measure space $(X,\mathcal{B},\mu)$, where $X=[0,6]$,
$\mathcal{B}$ is the Borel $\sigma$-field on $X$, and
$\mu(dx)=\left(\mu_{1},\mu_{2},\mu_{3}\right)(dx)=\left(f_{1}(x),f_{2}(x),f_{3}(x)\right)dx$,
where\[
\begin{array}{ccc}
f_{1}(x)=
\begin{cases}
30 & x\in[0,1),\\
40 & x\in[1,2),\\
10 & x\in[2,4),\\
15 & x\in[4,5),\\
5 & x\in[5,6];
\end{cases}
&
f_{2}(x)=
\begin{cases}
40 & x\in[0,1),\\
10 & x\in[1,2),\\
20 & x\in[2,4),\\
10 & x\in[4,5),\\
30 & x\in[5,6];
\end{cases}
&
f_{3}(x)=
\begin{cases}
10 & x\in[0,1),\\
20 & x\in[1,3),\\
30 & x\in[3,4),\\
20 & x\in[4,5),\\
25 & x\in[5,6].
\end{cases}
\end{array}
\]
These density functions are plotted in Fig.~\ref{fig:densities}.
\begin{figure}[h]
\begin{center}
\includegraphics[scale=0.39]{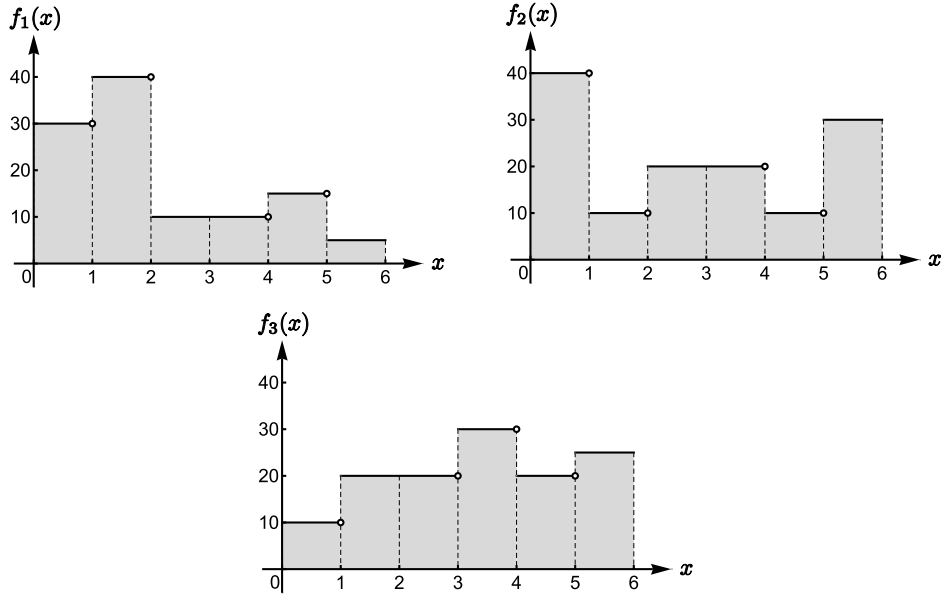} \label{fig:densities}
\caption{Density functions of the vector measure in
Example~\ref{ex:RinQ}.}
\end{center}
\end{figure}
Note that $\mu(X)=(110,130,125)$ and
\begin{equation}\label{e:rmpat} R_{\mu}(X)
=\left\{\sum\limits_{i=1}^6 \alpha_i p^i:\,\alpha_i \in [0,1], \
i=1,\dots,6\right\} \end{equation} is a zonotope, where $p^1 =
\mu([0,1)) = (30,40,10)$, $p^2 = \mu([1,2)) = (40,10,20)$, $p^3 =
\mu([2,3)) = (10,20,20)$, $p^4 = \mu([3,4)) = (10,20,30)$, $p^5 =
\mu([4,5)) = (15,10,20)$, $p^6 = \mu([5,6)) = (5,30,25)$.

 Let $p=p^1 + p^2 + p^3 =(80,70,50)$. Observe that $p$ is an extreme point
of $R_{\mu}(X)$. Indeed, consider the vector $d=\left( \frac{7}{5}, 1, -\frac{8}{5}  \right)$ and the linear function $l_d(\alpha)$
defined for all $\alpha\in {\mathbb R}^6$ by the scalar product
\begin{eqnarray*}
l_d(\alpha) &=& d \cdot \left( \sum\limits_{i=1}^6 \alpha_i p^i \right) = \sum\limits_{i=1}^6  \alpha_i (d \cdot p^i) \\
&=&  66 \alpha_1 + 34
\alpha_2 + 2 \alpha_3 -14 \alpha_4 - \alpha_5 - 3 \alpha_6.
\end{eqnarray*}

On the set $R_\mu(X)$, explicitly presented in (\ref{e:rmpat}),
this function achieves a unique maximum at the point
$\alpha^*=(1,1,1,0,0,0)$ with $l_d(\alpha^*)=66 + 34 + 2=102$. In
addition, $\sum_{i=1}^6\alpha^*_ip^i=p$.  So, $d\cdot r-102\le 0$
for all $r\in R_\mu(X)$ and the equality holds if and only if
$r=p.$ Thus, $d\cdot r-102=0$ is a supporting hyperplane of the
convex polytope $ R_\mu(X)$, and the intersection of the polytope
and hyperplane consists of the single point $p.$  This implies
that $p$ is an extreme point of $ R_\mu(X)$.

According to the definition of $R_\mu(X)$, for $p \in R_\mu(X)$
there exists a measurable subset $Z\in\mathcal{F}$ such that
$\mu(Z)=p$ and, according to  \cite[Theorem III]{Liapounoff:1940}
described in Section \ref{sec:intro}, since $p$ is extreme, such
$Z$ is unique up to null sets.  In particular, $p=\mu(Z)$ for
$Z=[0,3]$. Thus,
\[
R_{\mu}^{p}(X)=R_{\mu}(Z)=\left\{\sum_{i=1}^3 \alpha_ip^i:\,
\alpha_i\in [0,1],\ i=1,2,3 \right\}.
\]
Choose $q=(56,29,31)$ and observe that $q\notin R_{\mu}^{p}(X)$.
Indeed, $q\in R_{\mu}^{p}(X)$ if and only if there exist $\alpha_1, \alpha_2, \alpha_3 \in [0,1]$, such that $\sum_{i=1}^3
\alpha_i p^i=q$, which is equivalent to
\begin{eqnarray}\label{eq:linex}
\alpha_1(30,40,10)+\alpha_2(40,10,20)+\alpha_3(10,20,20)
=(56,29,31),
\end{eqnarray}
but the only solution to the linear system of equations
(\ref{eq:linex})  is
\[
\alpha_1=\frac{3}{10},\;\alpha_2=\frac{11}{10},\;\alpha_3=\frac{3}{10},\]
where $\alpha_2\notin[0,1]$.

On the other hand, $q\in Q_\mu^p(X),$  because: (i) $q\in
R_{\mu}(X)$, and (ii) $q\in R_{\mu}(X)-\left\{ \mu(X)-p\right\} $.
Indeed,  (i) holds since, for
$Z_{1}=\left[0,\frac{42}{115}\right)\cup\left[1,1\frac{229}{230}\right)\cup
\left[2,2\frac{33}{460}\right) \cup \left[4,4\frac{3}{10}\right)$,
\begin{eqnarray*}
\mu\left(Z_{1}\right) & = & \frac{42}{115}\times(30,40,10) + \frac{229}{230}\times(40,10,20)\\
 & + & \frac{33}{460}\times(10,20,20) +\frac{3}{10}\times(15,10,20)\\
 & = & (56,29,31)=q.
\end{eqnarray*}
Notice that (ii) is equivalent to $q+\mu(X)-p\in R_{\mu}(X)$,
where
$q+\mu(X)-p=(56,29,31)+(110,130,125)-(80,70,50)=(86,89,106)$. Let
$Z_{2}=\left[0,\frac{15}{46}\right)
\cup \left[1,1\frac{45}{46}\right)
\cup \left[2,2\frac{209}{230}\right)
\cup \left[3,5\right)
\cup \left[5,5\frac{3}{5}\right)$. Then
\begin{eqnarray*}
\mu\left(Z_{2}\right) & = & \frac{15}{46}\times(30,40,10)+\frac{45}{46}\times(40,10,20)+\frac{209}{230}\times(10,20,20)\\
 & + & 1\times(10,20,30)+1\times(15,10,20)+\frac{3}{5}\times(5,30,25) \\
 & = & (86, 89, 106) = q + \mu(X) - p.
\end{eqnarray*}
 Thus (ii) holds too, and $R_\mu^p(X)\ne Q_\mu^p(X)$. $\hfill \square$
\end{example}


The following example demonstrates that the necessary condition
(\ref{eq:necsuf}) for the existence of a measurable partition
$\{X^a:\, a\in A\}$ with $\mu(X^a)=p^a,$ $a\in A$,  is not
sufficient for an atomless measure $\mu$ when $m>2$ .  In this
example, $A$ consists of three points.  According to 
\cite[Theorem 2.5]{Dai:2010}, this condition is necessary and
sufficient when $m=2$, $A$ is countable, and $\mu$ is atomless. If
$A$ consists of two points, say $a$ and $b$, and $p^a\in
R_\mu(X)$, $p^b=\mu(X)-p^a$, then the partition $\{X^a,X^b\}$
always exists with $X^a$ selected as any $X^a\in{\mathcal F}$
satisfying $\mu(X^a)=p^a$ and with $X^b=X\setminus X^a.$

\begin{example} \label{ex:NoPart}
Consider the measure space $(X,\mathcal{B},\mu)$ defined in
Example \ref{ex:RinQ}. Let $p^{1}=(56,29,31)$, $p^{2}=(24,41,19)$,
$p^{3}=(30,60,75)$, and $A=\{1,2,3\}$. Then
$p^{1}+p^{2}+p^{3}=\mu(X)$. We further observe that: (i) $p^{1}$
is  the vector $q$ from Example \ref{ex:RinQ}, so $p^{1}\in
R_{\mu}(X)$ and therefore $p^2 + p^3 = \mu(X) - p^1 \in R_\mu (X)$; (ii) $p^{1}+p^{3}$ is the vector $q+\mu(X)-p$ from
Example \ref{ex:RinQ}, so $p^{1}+p^{3}\in R_{\mu}(X)$ and
therefore $p^2=\mu(X)-p^1-p^3\in R_\mu(X)$; (iii) $p^{1}+p^{2}$ is
the vector $p$ from Example \ref{ex:RinQ}, so $p^{1}+p^{2}\in
R_{\mu}(X)$ and therefore
$p^3=\mu(X)-p^1-p^2\in R_\mu(X)$. 
Thus, the vectors $p^{a}$, $a\in A$, satisfy (\ref{eq:necsuf}).

If there exists a partition $\{X^{a}\in\mathcal{B}:a\in A\}$ of
$X$ with $\mu(X^{a})=p^{a}$ for all $a\in A$, let $Y=X^{1}\cup
X^{2}$. Since $X^1\cap X^2=\emptyset$, $\mu(X^1)=p^1=q,$ and
$\mu(Y)=p^1+p^2=p, $ then $q\in R_\mu^p(X).$ According to
Example~\ref{ex:RinQ}, $q\notin R_\mu^p(X).$ This contradiction
implies that a partition $\{X^{a}\in\mathcal{B}:a\in A\}$ of $X$,
with $\mu(X^{a})=p^{a}$ for all $a\in A$, does not exist.$\hfill
\square$
\end{example}

In conclusion, we provide two simple examples  showing that, if
$\mu$ is not atomless, then even for $m=1$ (and, therefore, for
any natural number $m$) a maximal subset (\ref{eq:2D}) may not
exist and equation (\ref{eq:2DA}) may not hold.

\begin{example}\label{ex3}
Consider the measure space $(X,2^X,\mu)$, where $X=\{1,2,3,4\}$
and
\[
\mu(\{1\})=0.1, \; \mu(\{2\})=0.4, \; \mu(\{3\})=0.2, \;
\mu(\{3\})=0.3.
\]
Let $p=0.5$. Then $\mathcal{S}^{p}_\mu = \{ \{1,2\}, \{3,4\} \}$.
In other words, the only subsets that have the measure $0.5$ are
$Z^1=\{1,2\}$ and $Z^2=\{3,4\}$. However, $R_\mu(Z^1)$ is not a
subset of $R_\mu(Z^2)$ and vice versa. Therefore,  a maximal
subset does not exist for $p=0.5$.$\hfill \square$
\end{example}

\begin{example}\label{ex4}
Consider the measure space $(X,2^X,\mu)$, where $X=\{1,2,3\}$ and
\[
\mu(\{1\})=0.1, \; \mu(\{2\})=0.55, \; \mu(\{3\})=0.35.
\]
The range of $\mu$ on $X$ is
$R_\mu(X)=\{0,0.1,0.35,0.45,0.55,0.65,0.9,1\}$. Let $p=0.55$. Then
$Q^p_{\mu} = \{ 0, 0.1, 0.45, 0.55\}$ and $R^p_{\mu} = \{ 0.55
\}$. Thus $R^p_{\mu} \subset Q^p_{\mu}$, but $R^p_{\mu} \neq
Q^p_{\mu}$.$\hfill \square$
\end{example}

\end{document}